\documentclass[12pt,dvips]{article}
\usepackage[dvips]{graphicx}
\usepackage{color}
\usepackage{subfigure}
\usepackage{rotating}
\usepackage{amsmath}
\usepackage{amsthm}
\usepackage{amsfonts}

\newtheorem{theorem}{Theorem}[section]
\newtheorem{definition}[theorem]{Definition}

\makeatletter
\def\@biblabel#1{}
\makeatother
\begin{document}

%

\begin{center} {\Large{\bf The linear topology associated with weak convergence of probability measures}}
\bigskip


\bigskip

Liang Hong
\footnote{Liang Hong is an Assistant Professor in the Department of Mathematics, Robert Morris University, 6001 University Boulevard, Moon Township, PA 15108, USA. Tel.: (412) 397-4024. Email address: hong@rmu.edu.}

\end{center}

\vskip 10pt
\bigskip

\centerline{\noindent {\bf Abstract}} \noindent
This expository note aims at illustrating weak convergence of probability measures from a broader view than a previously published paper. Though the results are standard for functional analysts, this approach is rarely known by statisticians and our presentation gives an alternative view than most standard probability textbooks. In particular, this functional approach clarifies the underlying topological structure of weak convergence. We hope this short note is helpful for those who are interested in weak convergence as well as instructors of measure theoretic probability.

\bigskip

\noindent
{\it MSC 2010 Classification:} Primary 00-01; Secondary 60F05, 60A10. \\
\bigskip

\noindent {\it Keywords:} Weak convergence; probability measure; topological vector space; dual pair\\

\vskip 15 pt

\newpage
\section{Introduction}
Weak convergence of probability measures is often defined as follows (Billingsley 1999 or Parthasarathy 1967).
\begin{definition}\label{definition1.1}
Let $X$ be a metrizable space. A sequence $\{P_n\}$ of probability measures on $X$ is said to \textbf{converge weakly} to a probability measure $P$ if
\begin{equation*}
\lim_{n\rightarrow\infty}\int_X f(x) d P_n(dx)= \int_X f(x) P(dx)
\end{equation*}
for every bounded continuous function $f$ on $X$.
\end{definition}

It is natural to ask the following questions:
\begin{enumerate}
  \item [1.]Why the type of convergence defined in Definition \ref{definition1.1} is called ``weak convergence''?
  \item [2.]How weak it is?
  \item [3.]Is there any connection between weak convergence of probability measures and weak convergence in functional analysis?
  \item [4.]Why we do not use a metric to describe it?
\end{enumerate}
Answers to these questions are not obvious from Definition \ref{definition1.1}. Varadarajan (1958) gives answers to Questions 3 and 4 under the assumption that $X$ is either compact or separable. In addition, Varadarajan (1958) derives several important properties of the space of probability measures.
\begin{theorem}[Varadarajan 1958]\label{theorem1.3}
Suppose $X$ is a metrizable space and $\mathcal{P}(X)$ is the space of probability measures on $X$. Then the following statements hold.
\begin{enumerate}
  \item [(1)]$X$ is separable if and only if $\mathcal{P}(X)$ is separable.
  \item [(2)]$X$ is compact if and only if $\mathcal{P}(X)$ is compact.
  \item [(3)]$X$ is Polish if and only if $\mathcal{P}(X)$ is Polish.
\end{enumerate}
\end{theorem}
\noindent\textbf{Remark.} To be precise, here the topology on $\mathcal{P}(X)$ is the weak* topology $\sigma (\mathcal{P}(X), C_b(X))$ we shall define in Section 2.\\

In the next section, we answer Questions 1-4 using the theory of topological vector spaces. This approach has certain advantages: (1) it clearly illustrates the underlying topological structure of weak convergence; (2) it may lead to simple and clean proofs; (3) it allows one to work under minimal hypotheses. As a result, our framework is more general than that of Varadarajan (1958) since we only require $X$ to be a metrizable space, that is, $X$ need not be compact or separable.

\section{Weak convergence of probability measures: a topological vector space point of view}
Varadarajan (1958) explains weak convergence of probability measures from the Banach space point of view. In this section, we attempt to illustrate weak convergence from the topological vector space point of view. Since a topological vector space need not be a normed space, our setup is more general than that of Varadarajan (1958). In particular, we do not assume that $X$ is compact or separable. We will mainly follow Aliprantis and Border (2006) and Bourbaki (1987). At this point, we recommend that readers consult II. 40-43 of Bourbaki (1987) before proceeding further. Unless otherwise stated, $X$ denotes a metrizable space, $\mathcal{B}_X$ denotes the Borel $\sigma$-algebra on $X$, $\Psi(X)$ denotes the space of all finite signed measures on $(X, \mathcal{B}_X)$, $\mathcal{P}(X)$ denotes the space of all probability measures on $(X, \mathcal{B}_X)$, $C_b(X)$ denotes the space of all bounded continuous functions on $X$, and $C_u(X)$ denotes the space of all bounded $d$-uniformly continuous functions on $X$, where $d$ is a compatible metric. Define a map from $(C_b(X), \Psi(X))$ to $R$ by
\begin{equation}
(f, \mu)\mapsto \langle f, \mu\rangle\equiv\int_X f d\mu. \nonumber
\end{equation}
It is evident that the above map is linear in each variable separately, i.e., it is a bilinear form. Moreover, the following theorem shows that the space $C_u(X)$ and $\Psi(X)$ each separate points of the other. Since each bounded continuous function on $X$ can be approximated pointwisely by a sequence of bounded $d$-Lipschitz continuous functions (Corollary 3.15, Aliprantis and Border 2006), $C_b(X)$ and $\Psi(X)$ each separate points of the other.

\begin{theorem}[Varadrajan 1958]\label{theorem2.1}
For any two finite signed measures $\mu_1$ and $\mu_2$ on a metrizable space $X$, the following two conditions are equivalent:
\begin{enumerate}
  \item [(1)]$\mu_1=\mu_2$.
  \item [(2)]$\int_X f d\mu_1=\int_X f d\mu_2$ for every $f\in C_u(X)$.
\end{enumerate}
\end{theorem}

Therefore, $(C_b(X), \Psi(X))$ is a dual pair and $f\mapsto \langle f, \mu\rangle=\int_X fd\mu$ defines a functional from $\Psi(X)$ to $R$ (II. 40-41, Bourbaki 1987). Hence, we may identify $\Psi(X)$ as a subspace of $R^{C_b(X)}$ (equipped with the product topology) and obtain a Hausforff locally convex topology on $\Psi(X)$ (Proposition 2, II.43,  Bourbaki 1987), that is, the \textbf{weak* topology} $\sigma(\Psi (X), C_b(X))$. Since $\mathcal{P}(X)\subset \Psi (X)$, the weak* topology $\sigma(\Psi (X), C_b(X))$ induces a topology on $\mathcal{P}(X)$. This induced topology is often denoted by $\sigma(\mathcal{P}(X), C_b(X))$. In the literature of probability and statistics, $\sigma(\mathcal{P}(X), C_b(X))$ is often called the \textbf{weak topology} or \textbf{topology of convergence in distribution};  but precisely it is the induced weak* topology. From II. 42 of Bourbaki (1987), we know that a neighborhood base at $Q\in \mathcal{P}(X)$ is given by all sets of the form
\begin{equation}\label{2.2}
\left\{ P\in\mathcal{P}(X): \bigg|\int f_kd P-\int f_k dQ\bigg|\leq \epsilon, k=1, 2, ..., n\right\},
\end{equation}
where $n$ is a positive integer, $f_1, ..., f_n\in C_b(X)$, and $\epsilon>0$. From now on, $\mathcal{P}(X)$ will always be equipped with this topology unless stated otherwise. Indeed, we have given answers to Questions 1 and 3 in Section 1. Definition 2 on II. 42 of Bourbaki (1987) implies the answer to Question 2: weak convergence of probability measures corresponds to the weakest topology on $\mathcal{P}(X)$ that makes all the maps $f \mapsto \int f d P, \ P\in\mathcal{P}(X)$ continuous. Next, we answer Question 4.

In general, a topological vector space need not be metrizable. For a topological vector space, the metrizability condition is given by the following theorem (I. 16, Bourbaki 1987).

\begin{theorem}\label{theorem2.2}
A Hausdorff topological vector space is metrizable if and only if it has a countable neighborhood base at zero. \end{theorem}

This theorem might be a little bit too general since we have noticed that the induced weak* topology $\sigma(\mathcal{P}(X), C_b(X))$ is Hausdorff and locally convex. Indeed, the metrizability of a Hausdorff locally convex space can also be characterized (Corollary, II. 24, Bourbaki 1987).

\begin{theorem}\label{theorem2.1}
A Hausdorff locally convex space $(X, \mathcal{T})$ is metrizable if and only if $\mathcal{T}$  is generated by a sequence $\{q_n\}$ of seminorms; in this case, $\mathcal{T}$ is generated by the metric $d$ given by
\begin{equation*}
d(x, y)=\sum_{n=1}^{\infty}\frac{q_n(x-y)}{2^n[1+q_n(x-y)]}.
\end{equation*}
\end{theorem}
Clearly, the sufficient condition in Theorem \ref{theorem2.1} need not hold for the weak* topology $\sigma(\mathcal{P}(X), C_b(X))$. Therefore, in general we can only use neighborhood bases or subbases to describe weak topology. This answers Question 4 in Section 1.

Following the same argument in Section 2, the weak* topology $\sigma(\Psi (X), C_u(X))$ is well-defined; it induces a relative topology on $\mathcal{P}(X)$, which is often denoted as $\sigma(\mathcal{P}(X), C_u(X))$. Since $C_u(X)$ is a proper subset of $C_b(X)$ unless $X$ is compact, $\sigma(\Psi (X), C_b(X))$ is generally stronger than $\sigma(\Psi(X), C_u(X))$ (Corollary 3, II. 43, Bourbakai 1987). Surprisingly, the following theorem (Theorem 15.2, Aliprantis and Border 2006) shows that they induce the same topology on $\mathcal{P}(X)$.

\begin{theorem}\label{theorem3.1}
Let $X$ be a metrizable space and $d$ be a compatible metric. Then $\sigma(\mathcal{P}(X), C_b(X))=\sigma(\mathcal{P}(X), C_u(X))$.
\end{theorem}

Finally, we remark that the $\sigma(\Psi(X), C_b(X))$ admits a uniformity $\mathcal{D}$ on $\Psi(X)$ since a topological vector space is a commutative topological group. As a result, $\mathcal{P}(X)$ carries the corresponding relative uniformity $\mathcal{U}$. It is clear that a base for $\mathcal{D}$ consists of all the sets of the form
\begin{equation*}
\{(P_1, P_2)\mid P_1-P_2\in U\},
\end{equation*}
where $U$ is a neighborhood at zero in the space $(\Psi(X), \sigma(\Psi(X), C_b(X)))$ and that $\sigma(\mathcal{P}(X), C_b(X))$ is the topology on the uniform space $(\mathcal{P}(X),  \mathcal{U})$. If $\{P_n\}$ is a sequence  in $\mathcal{P}(X)$ and $\lim_n\int fdP_n$ exists for every $f\in C_b(X)$, then $\{P_n\}$ is $\mathcal{U}$-Cauchy. However, there is no guarantee that $\{P_n\}$ converges in $\sigma(\mathcal{P}(X), C_b(X))$ even if $\mathcal{P}(X)$ is Polish; because different uniformities may generate the same topology.

\newpage
\section*{Acknowledgments}
\noindent The author benefited from some discussions with Patrik J. Fitzsimmons, Lutz Mattner and Pietro Rigo on this topic. Thanks are also due to Christoph Lendenmann who read an earlier draft and caught some typos and obscurities.

\end{document}